\documentclass[a4paper,12pt,leqno]{amsart}

\usepackage{amsmath,amssymb,eucal,amsthm}
\usepackage{color}
\usepackage{bm}
\usepackage{graphicx}
\usepackage{mathrsfs}

\usepackage{latexsym}
\numberwithin{equation}{section}

\newtheorem{theorem}{Theorem}[section]
\newtheorem{lemma}[theorem]{Lemma}

\newtheorem{corollary}[theorem]{Corollary}

\theoremstyle{definition}
\newtheorem{definition}[theorem]{Definition}
\newtheorem{example}[theorem]{Example}

\newtheorem{remark}[theorem]{Remark}

\numberwithin{equation}{section}

\allowdisplaybreaks

\begin{document}
\hfill\texttt{\jobname.tex}\qquad\today

\title[Laurent series expansions of multiple zeta-functions]
{Laurent series expansions of multiple zeta-functions of
Euler-Zagier type at integer points}

\author{Kohji Matsumoto}

\address{K. Matsumoto, Graduate School of Mathematics Nagoya University, Chikusa-ku, 
Nagoya 464-8602, Japan, Email:{\rm kohjimat@math.nagoya-u.ac.jp}}

\author{Tomokazu Onozuka}

\address{T. Onozuka, Graduate School of Mathematics Nagoya University, Chikusa-ku,      
Nagoya 464-8602, Japan, Email:{\rm m11022v@math.nagoya-u.ac.jp}}

\footnote[0]{This research was partially 
supported by Grants-in-Aid for Scientific Research, Grant Numbers 25287002 (for the
first-named author) and 13J00312 (for the second-named author), JSPS.}

\author{Isao Wakabayashi}

\address{I. Wakabayashi, 3-19-8 Kugayama, Suginami-ku, Tokyo 168-0082, Japan, 
Email:{\rm wakaba.isao@gmail.com}}

\date{}

\maketitle

\begin{abstract}
We give explicit expressions (or at least an algorithm of obtaining such expressions)
of the coefficients of the Laurent series
expansions of the Euler-Zagier multiple zeta-functions at any integer points.
The main tools are the Mellin-Barnes integral formula and the harmonic product
formulas. The Mellin-Barnes integral formula is used in the induction process
on the number of variables, and the harmonic product formula is used to show that 
the Laurent series expansion outside the domain of convergence can be obtained from that 
inside the domain of convergence.
\end{abstract}

\section{Introduction and the statement of main results}\label{sec-1}

The Euler-Zagier $r$-ple zeta-function is defined by
(\cite{Hof92}, \cite{Zag94})
\begin{align}\label{1-1}
\zeta_r(\mathbf{s})=\sum_{n_1=1}^{\infty}\sum_{n_2=1}^{\infty}
\cdots\sum_{n_r=1}^{\infty}
n_1^{-s_1}(n_1+n_2)^{-s_2}\cdots(n_1+\cdots+n_r)^{-s_r}
\end{align}
(where $\mathbf{s}=(s_1,\ldots,s_r)\in{\mathbb C}^r$)
in the domain of its absolute convergence, which is
\begin{align}\label{1-2}
\mathcal{D}_r=\{\mathbf{s}\in\mathbb{C}^r\;|\;
\Re(s(j,r))>r-j+1\;(1\leq j\leq r)\},
\end{align}
where $s(j,r)=s_{j}+s_{j+1}+\cdots+s_r$
(\cite[Theorem 3]{Mat02}).
Special values $\zeta_r(\mathbf{m})$ ($\mathbf{m}=(m_1,\ldots,m_r)\in\mathbb{Z}^r$)
of \eqref{1-1} in this domain 
are called multiple zeta values (MZV),
and have been studied extensively.   

It is known that \eqref{1-1} can be continued meromorphically
to the whole space $\mathbb{C}^r$ (Akiyama, Egami and Tanigawa \cite{AET01},
Zhao \cite{Zha00}).     
Therefore the behavior of \eqref{1-1} around the points $\mathbf{m}\in\mathbb{Z}^r$
outside the domain $\mathcal{D}_r$ is also of great interest.   This direction of
research was initiated by Akiyama, Egami and Tanigawa \cite{AET01}, and then pursued
further by subsequent mathematicians
(Akiyama and Tanigawa \cite{AT01}, Komori \cite{Kom10}, Sasaki \cite{Sas09} \cite{Sas09b},
and the second-named author \cite{Ono13}). 

In this paper, to understand the behavior of \eqref{1-1} around the integer points more 
closely, we give Laurent series expansions for the Euler-Zagier $r$-ple zeta-function at integer 
points.  When $r=1$, the function $\zeta_1(s)$ is nothing but 
the Riemann zeta-function, and this function has, at $m\in\mathbb{Z}$, the following 
Taylor or Laurent series expansion:
\begin{align}\label{LSE1}
\zeta(s)&=\begin{cases}\displaystyle
\sum_{n=0}^\infty\frac{1}{n!}\zeta^{(n)}(m)(s-m)^n&(m\neq 1),\\
\displaystyle\frac{1}{s-1}+\sum_{n=0}^\infty\gamma_n(s-1)^n&(m=1),
\end{cases}
\end{align}
where $\gamma_n$ is the $n$-th Stieltjes constant (or generalized Euler constant)\footnote[2]{Some authors call $(-1)^nn!\gamma_n$ the $n$-th Stieltjes constant.} and $\zeta^{(n)}(s)$ is the $n$-th derivative of the Riemann zeta-function.  
The $0$-th Stieltjes constant $\gamma_0=\gamma$ is the well-known 
Euler constant, and Stieltjes constants have been studied by several authors; see \cite[p.164]{Ber85}. 
Note that for $m>1$, $\zeta^{(n)}(m)$ is given by the series $\sum_{k=1}^\infty (-\log k)^n/k^m$. (In this paper, we define $(-\log 1)^0=0^0=1$.) 

The same type of expansions holds in the multivariable case.
Generally, for a function $f(\mathbf{s})$ and $(n_1,\ldots,n_r)\in\mathbb{Z}_{\geq 0}$, 
we denote by $f^{(n_1,\ldots,n_r)}(\mathbf{s})$ 
the (partial) derivative 
$(\partial^{n_1}/\partial s_1^{n_1})\cdots(\partial^{n_r}/\partial s_r^{n_r})f(\mathbf{s})$.
When $\mathbf{m}\in\mathcal{D}_r$, obviously
\begin{align}\label{TaylorExp}
\zeta_r(\mathbf{s})=\sum_{n_1=0}^{\infty}
\cdots\sum_{n_r=0}^{\infty}\frac{1}{n_1!\cdots n_r!}\zeta_r^{(n_1,\ldots,n_r)}
(m_1,\ldots,m_r)(s_1-m_1)^{n_1}\cdots(s_r-m_r)^{n_r}\
\end{align}
when $\mathbf{s}$ is close to $\mathbf{m}$.

Also, we define the $r$-ple $(n_1,\ldots,n_r)$-th Stieltjes constant 
$\gamma_{(n_1,\ldots,n_r)}$ by the following formula which is valid when $\mathbf{s}$
is close to $(1,\ldots,1)$:
\begin{align}\label{LSE2}
\zeta_r(s_1,\ldots,s_r)&=\left(\prod_{k=2}^{r}\frac{1}{s_k+\cdots+s_r-(r-k+1)}\right)\left\{\frac{1}{s_1+\cdots+s_r-r}\right.\\
&\left.\qquad+\sum_{n_1=0}^{\infty}
\cdots\sum_{n_r=0}^{\infty}\gamma_{(n_1,\ldots,n_r)}(s_1-1)^{n_1}\cdots(s_r-1)^{n_r}\right\}.\notag
\end{align}
When $r=1$, (regarding the empty product as 1) we find that \eqref{LSE2} is reduced to
the second formula of \eqref{LSE1}, so $\gamma_{(n_1)}=\gamma_{n_1}$.

In Section \ref{sec-2}, we will prove the above \eqref{LSE2}, which we call 
the Laurent series expansion at $(s_1,\ldots,s_r)=(1,\ldots,1)$ of the Euler-Zagier 
$r$-ple zeta-function (see Lemma \ref{LEM1}).

\begin{remark}\label{rem-Laurent}
In the present paper we use the term "Laurent series expansion" in the above extended
way, that is, the linear factors in the denominator and in the numerator are not
necessarily the same.    The standard definition of Laurent series expansion is more
restricted; see \cite[Chapter I, Section 5, pp.88-90]{Fuk63} or 
\cite[Erstes Kapitel, Sektion 20, pp.68-70]{Osg29}.   
\end{remark}

\begin{remark}\label{rem-Stieltjes}
It is not clear what is the most suitable definition of multiple Stieltjes constants.
See Remark \ref{rem-Euler} at the end of Section \ref{sec-3}.
\end{remark}

The first main theorem in this paper is as follows.

\begin{theorem}\label{Main1}
Let $r,m_1,\ldots,m_r$ be positive integers, and put
$\mathbf{m}=(m_1,\ldots,m_r)$. 
Then the coefficients in the Laurent series of $\zeta_r(\mathbf{s})$ at $\mathbf{s}=\mathbf{m}$ can be given explicitly by using 
$\gamma_{(n_1,\ldots,n_k)}$ and $\zeta_k^{(l_1,\ldots,l_k)}(q_1,\ldots,q_k)$ where 
$1\leq k\leq r$, $(n_1,\ldots,n_k),(l_1,\ldots,l_k)\in(\mathbb{Z}_{\geq 0})^k$ and 
$(q_1,\ldots,q_k)\in\mathcal{D}_k\cap(\mathbb{Z}_{\geq 1})^k$.
\end{theorem}

Here we give some comments on the meaning of this theorem.

(i)
When $r=1$, this theorem is trivial, because the assertion is nothing but \eqref{LSE1}
for $m\geq 1$.
Similarly, in the general $r$-ple case
for $r\geq1$, the theorem is trivial by (\ref{TaylorExp}) when 
$\mathbf{m}\in\mathcal{D}_r\cap(\mathbb{Z}_{\geq 1})^r$,
and follows directly from (\ref{LSE2}) when $\mathbf{m}=(1,\ldots,1)$.
The main point of the theorem lies in the cases when 
$\mathbf{m}\in(\mathbb{Z}_{\geq 1})^r\setminus(\mathcal{D}_r\cup(1,\ldots,1))$, 
which happens only when $r\geq 2$.

(ii)
When 
$\mathbf{m}\in(\mathbb{Z}_{\geq 1})^r\setminus(\mathcal{D}_r\cup(1,\ldots,1))$,
we understand the meaning of the "Laurent series" at $\mathbf{s}=\mathbf{m}$ 
in the extended sense mentioned in Remark \ref{rem-Laurent}.
More strictly, when $m_j>1$ and $m_{j+1}=\cdots=m_r=1$ ($1\leq j\leq r-1$), then 
the Laurent series at
$\mathbf{s}=\mathbf{m}$ is of the form of a fraction, whose denominator is the product
of linear factors $(s(k,r)-(r-k+1))$ ($j+1\leq k\leq r$), and its numerator is a
Taylor series with respect to $(s_1-m_1),\ldots,(s_r-m_r)$.  
An example will be given in Example \ref{ex} at the end of Section \ref{sec-2}.  
The term "coefficients in the Laurent series" in the statement of the theorem 
means the coefficients in this Taylor series.
As can be seen in Example \ref{ex}, our proof of Theorem \ref{Main1} gives an
algorithm of obtaining the coefficients explicitly.
\bigskip

After mentioning some preparatory results in Section \ref{sec-1.5},
we will prove this theorem in Section \ref{sec-2}. 
Note that similar to the one variable case, $\zeta_k^{(l_1,\ldots,l_k)}(q_1,\ldots,q_k)$ also has the series expression
\begin{align*}
\sum_{n_1=1}^{\infty}\cdots\sum_{n_{k}=1}^{\infty} \frac{(-\log n_1)^{l_1}(-\log (n_1+n_2))^{l_2}\cdots(-\log(n_1+\cdots+n_k))^{l_k}}{n_1^{q_1}(n_1+n_2)^{q_2}\cdots(n_1+\cdots+n_k)^{q_k}}
\end{align*}
for $(l_1,\ldots,l_k)\in(\mathbb{Z}_{\geq 0})^k$ and $(q_1,\ldots,q_k)\in\mathcal{D}_k$.

In Section \ref{sec-3}, we will consider the coefficients in the Laurent series of $\zeta_r(\mathbf{s})$ under a certain restriction on the variables.  Under this restriction, these coefficients can be explicitly given by only using $\gamma_n$ and $\zeta_k^{(l_1,\ldots,l_k)}(q_1,\ldots,q_k)$. It implies that we can eliminate the role of 
$\gamma_{(n_1,\ldots,n_k)}$ from Theorem \ref{Main1} under this restriction. 

In Section \ref{sec-4}, we will consider the Laurent series expansions at 
$\mathbf{m}\in\mathbb{Z}^r\setminus(\mathbb{Z}_{\geq1})^r$. 
Since $\zeta_r(\mathbf{s})$ is singular on the hyperplanes
$$
s_r=1,\quad s_{r-1}+s_r=2,1,0,-2,-4,-6,\ldots,
$$
and
$$
s(j,r)=l\qquad(l\in\mathbb{Z},\;l\leq r-j+1)
$$
for $1\leq j\leq r-2$
(\cite[Theorem 1]{AET01}), the points
$\mathbf{m}\in \mathbb{Z}^r\setminus(\mathbb{Z}_{\geq1})^r$ 
are frequently on these singular hyperplanes, and in many cases are the points of indeterminacy.   All of the aforementioned previous studies on non-positive integer points
encountered this obstacle, and those studies discussed the limit values of
$\zeta_r(\mathbf{s})$ when $\mathbf{s}$ approaches $\mathbf{m}$ along various ways.    In particular, the second-named author \cite{Ono13} obtained a rather general result in which it allows a lot of flexibility how to approach the limit points, though he did not arrive at the Laurent series expansions. An example of this result will be mentioned just after Corollary \ref{prop3-1}.

We give Laurent series expansions at the points belonging to 
$\mathbb{Z}^r\setminus(\mathbb{Z}_{\geq1})^r$, 
which is the second main theorem in this paper.   Define
\begin{align}\label{def_M}
M_r({\mathbf m})=M_r(m_1,\ldots,m_r):=\max\{r-j-(m_j+\cdots+m_r)\;|\;1\leq j\leq r\}.
\end{align}
For $l\geq 2$ and $z_l\in\mathbb{C}$, let
\begin{align}\label{F-def}
F_{z_l}(\mathbf{s}):=\frac{\Gamma(s_l+z_l)}{\Gamma(s_l)}
\zeta_{l-1}(s_1,\ldots,s_{l-2},s_{l-1}+s_l+z_l)
\quad{\rm for}\;\;\mathbf{s}=(s_1,\ldots,s_l).
\end{align}
Then the result is the following
\begin{theorem}\label{Main3}
Let $\mathbf{m}\in\mathbb{Z}^r\setminus(\mathbb{Z}_{\geq1})^r$. Then the coefficients 
of the Laurent series of $\zeta_r(\mathbf{s})$ at $\mathbf{s}=\mathbf{m}$ can be given
explicitly in terms of the coefficients used in Theorem \ref{Main1}, 
$\zeta^{(n)}(m) (m\leq0,n\geq0)$, and integrals
\begin{align*}
\int_{(M_l(\mathbf{k})+1-\eta)}F_{z_l}^{(n_1,\ldots,n_l)}(\mathbf{k})\Gamma(-z_l)
\zeta(-z_l)dz_l
\end{align*}
where $2\leq l\leq r,\ (n_1,\ldots,n_l)\in(\mathbb{Z}_{\geq0})^l$, 
$\mathbf{k}\in\mathbb{Z}^l\setminus(\mathbb{Z}_{\geq1})^l$, $0<\eta<1$, and
the path of integration is the vertical line $\Re z_l=M_l(\mathbf{k})+1-\eta$.
\end{theorem}

A big difference from Theorem \ref{Main1} is that, here, the
coefficients may include some integrals. 
This theorem will be proved in the first half of Section \ref{sec-4}, and
in the second half of Section \ref{sec-4}, we will discuss some cases 
when we may ignore the contribution of integral terms.
\bigskip

The authors express their sincere gratitude to Professor Hidekazu Furusho,
Professor Yasushi Komori, Mr. Tomohiro Ikkai and Mr. Ryo Tanaka for valuable
comments and discussions.

\section{Preliminaries}\label{sec-1.5}

In the following sections,
we use in the induction process the following key formula proved by the first-named author (see \cite[(12.7)]{Mat03} or \cite[(4.4)]{Mat03b}):
\begin{align}\label{1-8}
\lefteqn{\zeta_r(s_1,\ldots,s_r)}\\
&=\frac{1}{s_r-1}\zeta_{r-1}(s_1,\ldots,s_{r-2},s_{r-1}+s_r-1)\notag\\
&\;+\sum_{k_r=0}^{M-1}\binom{-s_r}{k_r}\zeta_{r-1}
    (s_1,\ldots,s_{r-2},s_{r-1}+s_r+k_r)\zeta(-k_r)\notag\\
&\;+\frac{1}{\Gamma(s_r)}I(s_1,\ldots,s_r;M-\eta),\notag
\end{align}
where $r\geq 2$, $M$ is a positive integer, $\eta$ is a small positive number, and
\begin{align*}
\lefteqn{I(s_1,\ldots,s_r;\alpha)}\\
&=\frac{1}{2\pi i}\int_{(\alpha)}
\Gamma(s_r+z_r)\Gamma(-z_r)
\zeta_{r-1}(s_1,\ldots,s_{r-2},s_{r-1}+s_r+z_r)\zeta(-z_r)dz_r,
\end{align*}
whose path of integration is the vertical line $\Re z_r=\alpha$.
The formula (\ref{1-8}), which was proved by using the classical Mellin-Barnes integral formula, is valid in the region where the above integral is convergent.

For any points $\mathbf{m}\in\mathbb{Z}^r$, there exists a sufficiently large $M$ such that the integral of (\ref{1-8}) is analytic at $\mathbf{m}$ (see Lemma \ref{LEM2}). 
Hence this formula tells us the behavior of
$\zeta_r(s_1,\ldots,s_r)$ around the point $\mathbf{s}=\mathbf{m}$ from the information on the behavior of $\zeta_{r-1}$.

Another key formula which we use is the so-called harmonic product formula. 
This is used to show that the Laurent series expansion outside the domain of convergence can be obtained from that inside the domain of convergence.
The harmonic product formula is obtained by just decomposing the summation.
For example, the product of the Euler-Zagier double zeta-function and the Riemann 
zeta-function can be decomposed as follows:
\begin{align}\label{tripleA}
&\zeta_2(s_1,s_2)\zeta(s_3)\\
&=\sum_{n_1<n_2}\sum_{0<m}\frac{1}{n_1^{s_1}n_2^{s_2}m^{s_3}}\notag\\
&=\left(\sum_{n_1<n_2<m}+\sum_{n_1<n_2=m}+\sum_{n_1<m<n_2}+
\sum_{n_1=m<n_2}+\sum_{m<n_1<n_2}\right)\frac{1}{n_1^{s_1}n_2^{s_2}m^{s_3}}\notag\\
&=\zeta_3(s_1,s_2,s_3)+\zeta_2(s_1,s_2+s_3)+\zeta_3(s_1,s_3,s_2)+\zeta_2(s_1+s_3,s_2)+\zeta_3(s_3,s_1,s_2).\notag
\end{align}
Similarly we obtain
\begin{align}\label{tripleB}
&\zeta(s_1)\zeta_2(s_2,s_3)\\
&=\zeta_3(s_1,s_2,s_3)+\zeta_2(s_1+s_2,s_3)+\zeta_3(s_2,s_1,s_3)
+\zeta_2(s_2,s_1+s_3)+\zeta_3(s_2,s_3,s_1).\notag
\end{align}
The same method can be applied to the decomposition of more general product 
of two Euler-Zagier multiple zeta-functions.   For example, as a direct
generalization of \eqref{tripleA}, we have
\begin{align}\label{prf1}
&\zeta_{r-1}(s_1,\ldots,s_{r-1})\zeta(s_r)\\
\notag&=\zeta_r(s_1,\ldots,s_r)+\zeta_r(s_1,\ldots,s_{r-2},s_r,s_{r-1})+\cdots+\zeta_r(s_r,s_1,\ldots,s_{r-1})\\
\notag&+\zeta_{r-1}(s_1,\ldots,s_{r-2},s_{r-1}+s_r)+\zeta_{r-1}(s_1,\ldots,s_{r-3},s_{r-2}+s_r,s_{r-1})+\cdots\\
\notag&+\zeta_{r-1}(s_1+s_r,s_2,\ldots,s_{r-1}).
\end{align}
The most general form of the decomposition can be written as
\begin{align}\label{prf2}
&\zeta_{j}(s_1,\ldots,s_{j})\zeta_{r-j}(s_{j+1},\ldots,s_{r})\\
\notag&=\zeta_r(s_1,\ldots,s_r)+\zeta_r(s_1,\ldots,s_{j-1},s_{j+1},s_j,s_{j+2},s_{j+3},\ldots,s_r)+\cdots\\
\notag&+\zeta_r(s_{j+1},s_1,s_2,\ldots,s_{j-1},s_j,s_{j+2},\ldots,s_r)+\cdots\\
\notag&+\zeta_r(s_{j+1},\ldots,s_{r},s_1,\ldots,s_{j})\\
\notag&+(\text{the sum of $\zeta_l$ $(l<r)$})
\end{align}
for $1\leq j\leq r-1$.

\section{The Laurent series expansion at positive integer points}
\label{sec-2}

The main  aim of this section is to prove Theorem \ref{Main1}.
First, we determine the order of the pole of the Euler-Zagier multiple zeta-function.
For $\delta>0$, let 
$${\mathcal E}_j(\delta) 
=\{(s_1,\ldots,s_{j})\in\mathbb{C}^{j}\;|\;
\Re s_l>1-\delta \;(1\leq l\leq j)\}.$$
\begin{lemma}\label{LEM1}
For each $j\in\mathbb{Z}_{\geq 1}$, there exists a function 
$h_j(s_1,\ldots,s_{j+1})$, analytic in the region 
${\mathcal E}_{j+1}(\delta_j)$ with a sufficiently small positive constant 
$\delta_j$ {\rm (}depending only on $j${\rm )}, such that,
for any $r\in\mathbb{Z}_{\geq 2}$, the identity
\begin{align}\label{eq.lem1}
&\zeta_r(\mathbf{s})\\
\notag&=\frac{1}{(s_r-1)(s_{r-1}+s_r-2)\cdots(s_2+\cdots+s_r-(r-1))}\zeta(s_1+\cdots+s_r-(r-1))\\
\notag&+\frac{1}{(s_r-1)(s_{r-1}+s_r-2)\cdots(s_3+\cdots+s_r-(r-2))}h_1(s_1,s_2)\\
\notag&+\frac{1}{(s_r-1)(s_{r-1}+s_r-2)\cdots(s_4+\cdots+s_r-(r-3))}h_2(s_1,s_2,s_3)\\
\notag&+\cdots\\
\notag&+\frac{1}{s_r-1}h_{r-2}(s_1,\ldots,s_{r-1})\\
\notag&+h_{r-1}(s_1,\ldots,s_r),
\end{align}
holds in the region ${\mathcal E}_{r}(\delta_{r-1})$, especially at any point
$\mathbf{m}\in(\mathbb{Z}_{\geq1})^{r}$.
\end{lemma}

To prove Lemma \ref{LEM1}, we use the part (ii) of the following lemma.
\begin{lemma}\label{LEM2}
{\rm (i)}
For any $\mathbf{m}\in\mathbb{Z}^r$, if we choose $M\geq M_r({\mathbf m})+1$ 
{\rm (}where $M_r({\mathbf m})$ is defined in \eqref{def_M}{\rm )} and $0<\eta<1$, 
the integral $I(s_1,\ldots,s_r;M-\eta)$ is analytic at $\mathbf{m}$.

{\rm (ii)}
In particular, $I(s_1,\ldots,s_r;1-\eta)$ is analytic at any point
$\mathbf{m}\in(\mathbb{Z}_{\geq1})^{r}$.    In fact, 
$I(s_1,\ldots,s_r;1-\eta)$ is analytic in the region ${\mathcal E}_{r}(1/r)$.
\end{lemma}

($Proof$ $of$ $Lemma$ $\ref{LEM2}.$) 
(i) Let ${\mathcal F}_r(M,\eta)$ be the set of 
$(s_1,\ldots,s_r)\in{\mathbb C}^r$ satisfying
\begin{align}\label{2-1}
\Re(s_j+\cdots+s_r)>r-j-M+\eta \quad(1\leq j\leq r).
\end{align}
By \cite[Section 12]{Mat03}, the integral $I(s_1,\ldots,s_r;M-\eta)$ is analytic on 
${\mathcal F}_r(M,\eta)$. 
Consider the case $\mathbf{s}=\mathbf{m}$.
We see that $\mathbf{s}=\mathbf{m}\in{\mathcal F}_r(M,\eta)$ if
\begin{align}\label{2-1'}
M>r-j-(m_j+\cdots+m_r)+\eta \quad(1\leq j\leq r).
\end{align}
It is obvious that under the conditions $M\geq M_r({\mathbf m})+1$ and 
$0<\eta<1$, the inequality (\ref{2-1'}) holds. 
Hence $\mathbf{m}\in{\mathcal F}_r(M,\eta)$ holds for 
$M\geq M_r({\mathbf m})+1$. 

(ii) If $\mathbf{m}\in(\mathbb{Z}_{\geq1})^{r}$, then $M_r({\mathbf m})\leq-1$,
so we can choose $M=1$ in assertion (i).   This implies the first half of
assertion (ii).   When $\Re s_l>1-\delta$ ($1\leq l\leq r$), then 
\eqref{2-1} with $M=1$ is satisfied if $(r-j+1)(1-\delta)>r-j-1+\eta$, that is,
$2-\eta>(r-j+1)\delta$ for $1\leq j\leq r$.    This is valid if we choose
any $\delta$ satisfying
$0<\delta<(2-\eta)/r$, especially $\delta=1/r$.   The second half of (ii) 
hence follows.
{\hfill $\square$}\\

($Proof$ $of$ $Lemma$ $\ref{LEM1}.$)
We use the induction on $r$. First, we consider the case $r=2$. 
By (\ref{1-8}) with $M=1$, we have
\begin{align*}
\zeta_2(s_1,s_2)
&=\frac{1}{s_2-1}\zeta(s_1+s_2-1)-\frac{1}{2}\zeta(s_1+s_2)
+\frac{1}{\Gamma(s_2)}I(s_1,s_2;1-\eta),
\end{align*}
because $\zeta(0)=-1/2$.   
>From Lemma \ref{LEM2} (ii), the last term is analytic at $\mathbf{s}\in(\mathbb{Z}_{\geq1})^{2}$. Furthermore the second term is also analytic at $\mathbf{s}\in(\mathbb{Z}_{\geq1})^{2}$, since $\zeta(s)$ has only one pole at $s=1$.
Hence the case $r=2$ is done. (The sum of the second term and the third term gives $h_1(s_1,s_2)$, with $0<\delta_1\leq 1/2$.)

Now we assume that (\ref{eq.lem1}) holds for $r-1$, that is, we have
\begin{align}\label{eq.lem1'}
&\zeta_{r-1}(s_1,\ldots,s_{r-2},s_{r-1}')\\
\notag&=\frac{1}{(s_{r-1}'-1)\cdots(s_2+\cdots+s_{r-2}+s_{r-1}'-(r-2))}\zeta(s_1+\cdots+s_{r-1}'-(r-2))\\
\notag&+\frac{1}{(s_{r-1}'-1)\cdots(s_3+\cdots+s_{r-2}+s_{r-1}'-(r-3))}h_1(s_1,s_2)\\
\notag&+\cdots\\
\notag&+\frac{1}{s_{r-1}'-1}h_{r-3}(s_1,\ldots,s_{r-2})\\
\notag&+h_{r-2}(s_1,\ldots,s_{r-1}').
\end{align}
On the other hand, we use (\ref{1-8}) with $M=1$ to obtain
\begin{align}\label{1-8'}
&\zeta_r(s_1,\ldots,s_r)\\
&=\frac{1}{s_r-1}\zeta_{r-1}(s_1,\ldots,s_{r-2},s_{r-1}+s_r-1)-\frac{1}{2}\zeta_{r-1}(s_1,\ldots,s_{r-2},s_{r-1}+s_r)\notag\\
&+\frac{1}{\Gamma(s_r)}I(s_1,\ldots,s_r;1-\eta).\notag
\end{align}
Substituting (\ref{eq.lem1'}) to (\ref{1-8'}) with $s_{r-1}'=s_{r-1}+s_r-1$, we have
\begin{align*}
&\zeta_r(s_1,\ldots,s_r)\\
\notag&=\frac{1}{s_r-1}\biggl\{\frac{1}{(s_{r-1}+s_r-2)\cdots(s_2+\cdots+s_r-(r-1))}\zeta(s_1+\cdots+s_r-(r-1))\biggr.\\
\notag&+\frac{1}{(s_{r-1}+s_r-2)\cdots(s_3+\cdots+s_r-(r-2))}h_1(s_1,s_2)\\
\notag&+\cdots\\
\notag&+\frac{1}{(s_{r-1}+s_r-2)}h_{r-3}(s_1,\ldots,s_{r-2})\\
\notag&\biggl.+h_{r-2}(s_1,\ldots,s_{r-2},s_{r-1}+s_r-1)\biggr\}\\
\notag&-\frac{1}{2}\zeta_{r-1}(s_1,\ldots,s_{r-2},s_{r-1}+s_r)
+\frac{1}{\Gamma(s_r)}I(s_1,\ldots,s_r;1-\eta).\notag
\end{align*}
Since $(m_1,\ldots,m_{r-2},m_{r-1}+m_r)\in\mathcal{D}_{r-1}$ holds for 
$\mathbf{m}\in(\mathbb{Z}_{\geq 1})^r$, it is clear that $\zeta_{r-1}(s_1,\ldots,s_{r-2},s_{r-1}+s_r)$ is analytic at these points. Furthermore it follows from Lemma \ref{LEM2} (ii) that the last term is also analytic.  The term containing $h_{r-2}$ can be written as
$$
\frac{h_{r-2}(s_1,\ldots,s_{r-2},s_{r-1}+s_r-1)-h_{r-2}(s_1,\ldots,s_{r-2},s_{r-1})}{s_r-1}+\frac{h_{r-2}(s_1,\ldots,s_{r-1})}{s_r-1}.
$$
The first term here is analytic at $\mathbf{m}\in(\mathbb{Z}_{\geq 1})^r$.
(When $s_r=1$, this term is to be understood as the derivative of
$h_{r-2}$ with respect to the last variable.)
Hence putting 
\begin{align*}
&h_{r-1}(s_1,\ldots,s_r)\\
&:=\frac{h_{r-2}(s_1,\ldots,s_{r-2},s_{r-1}+s_r-1)-h_{r-2}(s_1,\ldots,s_{r-1})}{s_r-1}\\
&-\frac{1}{2}\zeta_{r-1}(s_1,\ldots,s_{r-2},s_{r-1}+s_r)\\
&+\frac{1}{\Gamma(s_r)}I(s_1,\ldots,s_r;1-\eta),
\end{align*}
which is analytic in ${\mathcal E}_{r}(\delta_{r-1})$ for a sufficiently small
$\delta_{r-1}>0$, 
we obtain Lemma \ref{LEM1}.
{\hfill $\square$}\\

($Proof$ $of$ $\eqref{LSE2}.$)
This is immediate from the formula \eqref{eq.lem1} of Lemma \ref{LEM1}, 
with expanding the Riemann zeta factor to
the Laurent series, and $h_j$ factors to the Taylor series.
{\hfill $\square$}\\

Next we prove Theorem \ref{Main1}.\\

($Proof$ $of$ $Theorem$ $\ref{Main1}.$) We use the induction on $r$. 
When $r=1$, it follows from (\ref{LSE1}) that Theorem \ref{Main1} holds. 
We assume that Theorem \ref{Main1} holds 
when the number of variables is $1,2,\ldots,r-1$, and we prove that Theorem \ref{Main1} 
holds for $r$.

Let $\mathbf{m}\in(\mathbb{Z}_{\geq 1})^r$. 
We define the case $(C_j)$ as
\begin{align}
(C_j)\ :\ \left\{
\begin{array}{ll}
m_r>1\quad{\rm if}\;\; j=r,\\
m_j>1,m_{j+1}=m_{j+2}=\cdots=m_r=1 \quad{\rm if}\;\; 1\leq j\leq r-1.
\end{array}
\right.
\end{align}
When $m_r>1$, since 
$\mathbf{m}\in\mathcal{D}_r$ holds, the series (\ref{1-1}) is absolutely convergent 
at $\mathbf{m}$. Therefore $\zeta_r(\mathbf{s})$ has the Taylor series expansion
\eqref{TaylorExp} at $\mathbf{m}$, so we are done in the case $(C_r)$.

Next we consider the case $(C_j)$, $1\leq j\leq r-1$. 
We use the (down-going) induction on $j$.
Let $1\leq j\leq r-1$, 
assume that Theorem \ref{Main1} holds in the cases $(C_{j+1}),\ldots,(C_{r})$, and we prove Theorem \ref{Main1} for $(C_j)$.
Our tool is the formula \eqref{prf2}.
The first term on the right-hand side of \eqref{prf2} is that we want to expand. 
Coefficients in the Laurent series of the left-hand side and of the terms 
of $\zeta_{l}$ ($l<r$) can be given explicitly by using $\gamma_{(n_1,\ldots,n_k)}$ 
and $\zeta_k^{(l_1,\ldots,l_k)}(q_1,\ldots,q_k)$ by the assumption of 
induction on $r$. 
Coefficients in the Laurent series of the terms of $\zeta_r$ except 
the first term on the right-hand side can be also given explicitly by 
using $\gamma_{(n_1,\ldots,n_k)}$ and $\zeta_k^{(l_1,\ldots,l_k)}(q_1,\ldots,q_k)$,
by the assumption of induction on $j$, because in these terms, $s_j$ is 
located at the $l$-th element of the $r$-tuple where $l>j$.  
Hence in this case, Theorem \ref{Main1} holds.


Therefore by induction, we obtain Theorem \ref{Main1} for 
$(C_j)$ $(j=1,\ldots,r-1)$. Finally, the only remaining case 
$\mathbf{m}=(1,\ldots,1)$ is implied by (\ref{LSE2}).
Thus we complete the proof of Theorem \ref{Main1}.
(We can also see the fact mentioned in comment (ii) just after the statement of 
Theorem \ref{Main1}, by analyzing the above proof a little more carefully.)
{\hfill $\square$}\\

\begin{example}\label{ex}
Here we explain the procedure given in the proof of Theorem \ref{Main1}, by describing
the case $r=3$, $\mathbf{m}=(2,1,1)$.    Assume $s_1$ is close to 2, and $s_2, s_3$ are 
close to 1.    We first use \eqref{tripleB}:
\begin{align}\label{ex1}
&\zeta_3(s_1,s_2,s_3)=\zeta(s_1)\zeta_2(s_2,s_3)-\zeta_2(s_1+s_2,s_3)\\
&\;\;-\zeta_3(s_2,s_1,s_3)
-\zeta_2(s_2,s_1+s_3)-\zeta_3(s_2,s_3,s_1).\notag
\end{align}
Since $s_1$ is close to 2 and $s_1+s_3$ is close to 3, the last two terms on the 
right-hand side are in the domain of absolute convergence, so can be expanded by
\eqref{TaylorExp}.   In particular, these are $O(1)$.
Next, we apply \eqref{tripleA} 
to the third term on the right-hand side:
\begin{align}\label{ex2}
&\zeta_3(s_2,s_1,s_3)=\zeta(s_2,s_1)\zeta(s_3)-\zeta_2(s_2,s_1+s_3)\\
&\;\;-\zeta_3(s_2,s_3,s_1)-\zeta_2(s_2+s_3,s_1)-\zeta_3(s_3,s_2,s_1).\notag
\end{align}
Since $s_1$ is close to 2, all terms but the first one on the right-hand side are
in the domain of absolute convergence, hence can be expanded by \eqref{TaylorExp} and 
$O(1)$.   The first term can be expanded by \eqref{LSE1} and \eqref{TaylorExp}, and can
be written as
\begin{align*}
\zeta(s_2,s_1)\zeta(s_3)=\zeta_2(s_2,s_1)\left(\frac{1}{s_3-1}+O(1)\right).
\end{align*}
Therefore
\begin{align}\label{ex3}
\zeta_3(s_2,s_1,s_3)=\frac{\zeta_2(s_2,s_1)}{s_3-1}+O(1).
\end{align}
To the second term on the right-hand side of \eqref{ex1}, we apply the simplest
harmonic product formula
\begin{align}\label{harmonic}
\zeta(s_1)\zeta(s_2)=\zeta_2(s_1,s_2)+\zeta_2(s_2,s_1)+\zeta(s_1+s_2)
\end{align}
to obtain
\begin{align*}
\zeta_2(s_1+s_2,s_3)=\zeta(s_1+s_2)\zeta(s_3)-\zeta_2(s_3,s_1+s_2)
-\zeta(s_1+s_2+s_3).
\end{align*}
The second and the third terms on the right-hand side are in the domain of absolute convergence, and so
\begin{align}\label{ex4}
\zeta_2(s_1+s_2,s_3)=\frac{\zeta(s_1+s_2)}{s_3-1}+O(1).
\end{align}
Finally, since $(s_2,s_3)$ is close to $(1,1)$, we use \eqref{LSE2} to obtain
\begin{align}\label{ex5}
\zeta(s_1)\zeta_2(s_2,s_3)
=\zeta(s_1)\frac{1}{s_3-1}\left(\frac{1}{s_2+s_3-2}+A(s_2,s_3)\right),
\end{align}
where
$$
A(s_2,s_3)=\sum_{n_1=0}^{\infty}\sum_{n_2=0}^{\infty}
\gamma_{(n_1,n_2)}(s_2-1)^{n_1}(s_3-1)^{n_2}.
$$
By the above argument it is clear that $\zeta_3(s_1,s_2,s_3)$ can be expanded to
the Laurent series around the point $(2,1,1)$, and especially
\begin{align}\label{ex6}
&\zeta_3(s_1,s_2,s_3)=\frac{\zeta(s_1)}{(s_3-1)(s_2+s_3-2)}\\
&\qquad+\frac{\zeta(s_1)A(s_2,s_3)-\zeta(s_1+s_2)-\zeta_2(s_2,s_1)}{s_3-1}+O(1),
\notag
\end{align}
where the terms on the numerators are all holomorphic there.
Expanding the numerators and the $O(1)$ term to the Taylor series, we 
obtain an example of comment (ii)
after the statement of Theorem \ref{Main1}.
\end{example}

\section{The Laurent series expansion at positive integer points under a
certain additional restriction}\label{sec-3}

In this section, we consider the Laurent series expansion of \eqref{1-1}
under a certain additional restriction on the variables.

\begin{definition}
Let $\mathbf{m}=(m_1,\ldots,m_r)\in (\mathbb{Z}_{\geq1})^r$. 
Let $h$ be the number of $m_j$ which is equal to 1, and denote those $m_j$s
by $m_{a_1},\ldots,m_{a_h}$.
By the {\it restricted Laurent series expansion} of $\zeta_r(\mathbf{s})$ 
at $\mathbf{m}$ we mean the Laurent series expansion of 
$\zeta_r(\mathbf{s})$ at $\mathbf{m}$ with the restriction $s_{a_1}=\cdots=s_{a_h}=s$.\end{definition}
If we add the above restriction, then the coefficients of the Laurent series can be given in the following simpler form than Theorem \ref{Main1}.

\begin{theorem}\label{Main2}
Let $r,m_1,\ldots,m_r$ be positive integers. 
The coefficients in the restricted Laurent series of $\zeta_r(\mathbf{s})$ at $\mathbf{s}=\mathbf{m}$ can be given explicitly by using $\gamma_n$ and $\zeta_k^{(l_1,\ldots,l_k)}(q_1,\ldots,q_k)$ for $n\geq0,\ 1\leq k\leq r$, $(l_1,\ldots,l_k)\in(\mathbb{Z}_{\geq 0})^k$ and $(q_1,\ldots,q_k)\in\mathcal{D}_k\cap(\mathbb{Z}_{\geq 1})^k$.
\end{theorem}

($Proof$ $of$ $Theorem$ $\ref{Main2}.$)
The proof of Theorem \ref{Main2} is quite similar to the proof of Theorem 
\ref{Main1}, so we omit the details.   The main difference is the last step 
of induction, i.e.  the case $\mathbf{m}=(1,\ldots,1)$. In this step, we used (\ref{LSE2}) in the proof of Theorem \ref{Main1}. However, in the proof of Theorem \ref{Main2}, we can not use (\ref{LSE2}), since its coefficients 
include $\gamma_{(n_1,\ldots,n_k)}$. Therefore, we have to use a different method.

The idea is to use (\ref{prf1}) once more. In this case, since 
$s_1=\cdots=s_r=s$, we can simplify (\ref{prf1}) as follows:
\begin{align*}
&\zeta_{r-1}(s,\ldots,s)\zeta(s)=r\zeta_r(s,\ldots,s)\\
&\qquad+\zeta_{r-1}(s,\ldots,s,2s)+\zeta_{r-1}(s,\ldots,s,2s,s)+\cdots+\zeta_{r-1}(2s,s,\ldots,s).
\end{align*}
The first term of the right-hand side is that we want to expand, and other terms can be expanded by assumption. Hence we obtain Theorem \ref{Main2}.
{\hfill $\square$}\\

>From (\ref{LSE2}) and Theorem \ref{Main2}, we can easily deduce the following corollary.

\begin{corollary}\label{cor}
For any non-negative integer $N$, the sum 
$$
\sum_{\substack{n_1+\cdots+n_r=N\\n_1,\ldots,n_r\geq0}} \gamma_{(n_1,\ldots,n_r)}
$$ 
can be written explicitly by using  $\gamma_n$ and 
$\zeta_k^{(l_1,\ldots,l_k)}(q_1,\ldots,q_k)$ for $n\geq0,\ 1\leq k\leq r$, 
$(l_1,\ldots,l_k)\in(\mathbb{Z}_{\geq 0})^k$ and 
$(q_1,\ldots,q_k)\in\mathcal{D}_k\cap(\mathbb{Z}_{\geq 1})^k$.
\end{corollary}

($Proof$ $of$ $Corollary$ $\ref{cor}.$)
By (\ref{LSE2}), we have
\begin{align*}
&\zeta_r(s,\ldots,s)\\
&=\frac{1}{r!}\frac{1}{(s-1)^r}+\frac{1}{(r-1)!}\frac{1}{(s-1)^{r-1}}\sum_{n_1=0}^{\infty}
\cdots\sum_{n_r=0}^{\infty}\gamma_{(n_1,\ldots,n_r)}(s-1)^{n_1+\cdots+n_r}\\
&=\frac{1}{r!}\frac{1}{(s-1)^r}+\frac{1}{(r-1)!}\frac{1}{(s-1)^{r-1}}\sum_{N=0}^{\infty}\sum_{\substack{n_1+\cdots+n_r=N\\n_1,\ldots,n_r\geq0}}\gamma_{(n_1,\ldots,n_r)}(s-1)^{N}.
\end{align*}
The assertion of the corollary follows from Theorem \ref{Main2}, because 
the above summation is the restricted Laurent series expansion at $(1,\ldots,1)$.
{\hfill $\square$}\\

\begin{remark}\label{rem-Euler}
In the present paper we define multiple Stieltjes constants by \eqref{LSE2},
but it is not clear whether this is the most suitable definition, or not.
Here we mention another definition of double Euler constant $\gamma_2(s_1)$
($\Re s_1>1$), due to \cite[Section 5]{MT15}.

In \cite{MT15}, $\gamma_2(s_1)$ is defined as a double analogue of the classical
infinite series definition of $\gamma$, but according to \cite[Proposition 5.1]{MT15},
\begin{align}\label{formula-Euler}
\gamma_2(s_1)=\lim_{s_2\to 1}\left(\zeta_2(s_1,s_2)-\frac{\zeta(s_1)}{s_2-1}
\right)
\end{align}
(and also $\gamma_2(s_1)=\zeta(s_1)\gamma-\zeta_2(1,s_1)-\zeta(s_1+1)$)
holds.    On the other hand, when $r=2$, \eqref{LSE2} implies
$$
\zeta_2(s_1,s_2)=\frac{1}{s_2-1}\left\{\frac{1}{s_1+s_2-2}+\sum_{n_1=0}^{\infty}
\sum_{n_2=0}^{\infty}\gamma_{(n_1,n_2)}(s_1-1)^{n_1}(s_2-1)^{n_2}\right\}
$$
around the point $(1,1)$.
When $|s_2-1|<|s_1-1|$, substituting the expansion
$$
\frac{1}{s_1+s_2-2}=\frac{1}{s_1-1}\sum_{n_2=0}^{\infty}\left(-\frac{s_2-1}{s_1-1}
\right)^{n_2}
$$
into the above, we have
\begin{align}\label{doubleexpansion}
\zeta_2(s_1,s_2)&=\frac{1}{s_2-1}\sum_{n_2=0}^{\infty}A_{n_2}(s_1)(s_2-1)^{n_2},
\end{align}
where
\begin{align}\label{def-A}
A_{n_2}(s_1)=\frac{(-1)^{n_2}}{(s_1-1)^{n_2+1}}+
\sum_{n_1=0}^{\infty}\gamma_{(n_1,n_2)}(s_1-1)^{n_1}.
\end{align}
>From \eqref{doubleexpansion} we see that $(s_2-1)\zeta_2(s_1,s_2)$ tends to
$A_0(s_1)$ as $s_2\to 1$, but from \eqref{formula-Euler} we also see that the
same limit tends to $\zeta(s_1)$.   That is, $A_0(s_1)=\zeta(s_1)$.
(In particular, from the case $n_2=0$ of \eqref{def-A} we find that 
$\gamma_{(n_1,0)}=\gamma_{n_1}$.)
Therefore \eqref{doubleexpansion} can be rewritten as
$$
\zeta_2(s_1,s_2)-\frac{\zeta(s_1)}{s_2-1}=\frac{1}{s_2-1}
\sum_{n_2=1}^{\infty}A_{n_2}(s_1)(s_2-1)^{n_2},
$$
so with \eqref{formula-Euler} we obtain $\gamma_2(s_1)=A_1(s_1)$.
\end{remark}

\section{The Laurent series expansion at other integer points}\label{sec-4}

In this section, we consider the Laurent series at $\mathbf{m}\in\mathbb{Z}^r\setminus(\mathbb{Z}_{\geq1})^r$. It is difficult to treat these points, since the harmonic product does not work well. Therefore we only use the key formula (\ref{1-8}). By Lemma \ref{LEM2} (i), the integral in (\ref{1-8}) is analytic at $\mathbf{m}$ 
when we choose $M=M_r(\mathbf{m})+1$. Therefore, when
$\mathbf{s}$ is close to $\mathbf{m}$, we have
\begin{align}\label{LSE3}
\lefteqn{\zeta_r(s_1,\ldots,s_r)}\\
&=\frac{1}{s_r-1}\zeta_{r-1}(s_1,\ldots,s_{r-2},s_{r-1}+s_r-1)\notag\\
&\;+\sum_{k_r=0}^{M_r(\mathbf{m})}\binom{-s_r}{k_r}\zeta_{r-1}
    (s_1,\ldots,s_{r-2},s_{r-1}+s_r+k_r)\zeta(-k_r)\notag\\
&\;+\frac{1}{\Gamma(s_r)}I(s_1,\ldots,s_r;M_r(\mathbf{m})+1-\eta).
\notag
\end{align}
Since the integral term $I(s_1,\ldots,s_r;M_r(\mathbf{m})+1-\eta)$ is 
analytic at $\mathbf{m}$, we have
\begin{align}\label{LSE4}
&\frac{1}{\Gamma(s_r)}I(s_1,\ldots,s_r;M_r(\mathbf{m})+1-\eta)
\notag\\
&=\sum_{n_1=0}^{\infty}\cdots\sum_{n_r=0}^{\infty} \frac{1}{2\pi in_1!\cdots n_r!}\int_{(M_r(\mathbf{m})+1-\eta)}F_{z_r}^{(n_1,\ldots,n_r)}(\mathbf{m})\Gamma(-z_r)
\zeta(-z_r)dz_r\\
\notag&\qquad\qquad\qquad\qquad\qquad\qquad\qquad\qquad\qquad\times(s_1-m_1)^{n_1}\cdots(s_r-m_r)^{n_r},
\end{align}
where $F_{z_l}(\mathbf{s})$ is defined by \eqref{F-def}. The above change of summation and integration is possible because $\Gamma(-z_r)$ decays rapidly as $|\Im z_r|\to\infty$.

The first aim of this section is to give a proof of Theorem \ref{Main3}.
\bigskip

($Proof$ $of$ $Theorem$ $\ref{Main3}.$)
We use induction for $r$. When $r=1$, the theorem follows from (\ref{LSE1}). We assume that the theorem holds when the number of variables is $1,2,\ldots,r-1$, 
and we prove that the theorem holds for $r$. 

In the case $r$, we prove that it follows from (\ref{LSE3}) and (\ref{LSE4}) that Theorem \ref{Main3} holds. From (\ref{LSE3}), $\zeta_r(s_1,\ldots,s_r)$ can be decomposed as a sum of three terms which contain $\zeta_{r-1}$. The last term can be expanded by (\ref{LSE4}), hence we only consider the first two terms. These terms consist of binomial coefficients and $\zeta_{r-1}(s_1,\ldots,s_{r-2},s_{r-1}+s_r+k)$ for integer $k$. Binomial coefficients are just polynomials in $s_r$ with rational coefficients. 
When $(m_1,\ldots,m_{r-2},m_{r-1}+m_r+k)\in(\mathbb{Z}_{\geq1})^{r-1}$, we use Theorem \ref{Main1}. Then $\zeta_{r-1}(s_1,\ldots,s_{r-2},s_{r-1}+s_r+k)$ can be expanded by using the coefficients of Theorem \ref{Main1}. When $(m_1,\ldots,m_{r-2},m_{r-1}+m_r+k)\in\mathbb{Z}^{r-1}\setminus(\mathbb{Z}_{\geq1})^{r-1}$, we use the assumption of the induction. Thus, noting the
following remark, we obtain Theorem \ref{Main3}.

The remark is as follows.  We expand 
$\zeta_{r-1}(s_1,\ldots,s_{r-2},s_{r-1}+s_r+k)$ at the point 
$(m_1,\ldots,m_{r-2},m_{r-1}+m_r+k)$.   Then the shape of the numerator part of the
expansion is
$$
\sum(\mbox{coefficient})(s_1-m_1)^{a_1}\cdots(s_{r-2}-m_{r-2})^{a_{r-2}}(s_{r-1}+s_r-m_{r-1}-m_r)^{a_{r-1}}.
$$
Therefore we need to modify the last factor by using the relation
$$
(s_{r-1}+s_r-m_{r-1}-m_r)^{a_{r-1}}=\sum_{n=0}^{a_{r-1}}\binom{a_{r-1}}{n}(s_{r-1}-m_{r-1})^n(s_r-m_r)^{a_{r-1}-n}.
$$
{\hfill $\square$}\\

In Theorem \ref{Main3}, the coefficients are not simple, since they include integrals whose behaviour is not obvious.   
However,
since the gamma function has poles at non-positive integers, we can estimate the last term of (\ref{LSE3}) by $O(|s_r-m_r|)$ for 
$m_r\in\mathbb{Z}_{\leq0}$.    Therefore in this case
\begin{align}
\lefteqn{\zeta_r(s_1,\ldots,s_r)}\notag\\
&=\frac{1}{s_r-1}\zeta_{r-1}(s_1,\ldots,s_{r-2},s_{r-1}+s_r-1)\notag\\
&\;+\sum_{k_r=0}^{M_r(\mathbf{m})}\binom{-s_r}{k_r}\zeta_{r-1}
    (s_1,\ldots,s_{r-2},s_{r-1}+s_r+k_r)\zeta(-k_r)\notag\\
&\;+O(|s_r-m_r|).\notag
\end{align}

Generally, the coefficients of the Laurent series expansion of $\zeta_{r-1}(s_1,\ldots,s_{r-2},s_{r-1}+s_r+k_r)$ also include integrals, and we cannot ignore the behavior of those terms. However, in some cases, such terms do not appear.   Hereafter we discuss some examples when we can, or cannot, obtain the limit value if we regard the term involving the integral $I(s_1,\ldots,s_r;M_r(\mathbf{m})+1-\eta)$ as an $O$-term.    For example, in the case $r=2$ and $m_2\leq0$, 
since $\zeta_1$ is the Riemann zeta function, we have
\begin{align}
\zeta_2(s_1,s_2)&=\frac{1}{s_2-1}\zeta(s_1+s_2-1)\notag\\
&\;+\sum_{k_2=0}^{M_2(\mathbf{m})}\binom{-s_2}{k_2}\zeta(s_1+s_2+k_2)\zeta(-k_2)\notag\\
&\;+O(|s_2-m_2|).\notag
\end{align}
This equation determines the limit value at the point $(m_1,m_2)\in\mathbb{Z}\times\mathbb{Z}_{\leq0}$.  That is, we obtain the following theorem.
\begin{theorem}\label{Main4}
For $(m_1,m_2)\in\mathbb{Z}\times\mathbb{Z}_{\leq0}$ and $\varepsilon_1,\varepsilon_2\in\mathbb{C}$, we have
\begin{align}\label{2-7}
\lefteqn{\zeta_2(m_1+\varepsilon_1,m_2+\varepsilon_2)}\\
&=\frac{1}{m_2-1+\varepsilon_2}\zeta(m_1+m_2-1+\varepsilon_1+\varepsilon_2)\notag\\
&\;+\sum_{k_2=0}^{M_2({\mathbf m})}\binom{-m_2-\varepsilon_2}{k_2}\zeta(-k_2)
\zeta(m_1+m_2+k_2+\varepsilon_1+\varepsilon_2)\notag\\
&\;+O(|\varepsilon_2|).\notag
\end{align}
\end{theorem}
In particular, the following formula holds.

\begin{corollary}\label{prop3-1}
When $m_1\leq 0$, $m_2\leq 0$, we have
\begin{align*}
\lefteqn{\zeta_2(m_1+\varepsilon_1,m_2+\varepsilon_2)}\\
&=\frac{1}{m_2-1}\zeta(m_1+m_2-1)\notag\\
&\;+\sum_{k_2=0}^{-m_2}
\frac{(-1)^{k_2}}{k_2!}m_2(m_2+1)\cdots(m_2+k_2-1)\zeta(-k_2)\zeta(m_1+m_2+k_2)
\notag\\
&\;+\frac{(-1)^{1-m_1-m_2}}{(1-m_1-m_2)!}\frac{(m_2+\varepsilon_2)                      
(m_2+1+\varepsilon_2)\cdots(-m_1+\varepsilon_2)}{\varepsilon_1+\varepsilon_2}           
\zeta(m_1+m_2-1)\notag\\                                                                
&\;+O(|\varepsilon_2|)+O(|\varepsilon_1+\varepsilon_2|).\notag
\end{align*}
\end{corollary}
As mentioned in the introduction, this kind of limit values was studied by several mathematicians (\cite{AET01}\cite{AT01}\cite{Kom10}\cite{Sas09}\cite{Sas09b}).
The second-named author also studied such limit values in \cite{Ono13}, under 
the conditions
\begin{align}\label{1-5}
\varepsilon_2\neq0\ ,\varepsilon_1+\varepsilon_2\neq 0
\end{align}
and
\begin{align}\label{1-6}
\left|\frac{\varepsilon_1}{\varepsilon_1+\varepsilon_2}\right|\ll 1,\ \left|\frac{\varepsilon_2}{\varepsilon_1+\varepsilon_2}\right|\ll 1,
\end{align}
which are weaker than the conditions assumed in the previous works.

An advantage of the above corollary is that even the condition 
\eqref{1-6} is not required.

($Proof$ $of$ $Corollary$ $\ref{prop3-1}.$)
When $m_1\leq 0$, $m_2\leq 0$, the case $m_1+m_2=2$ does not occur, so the first term of (\ref{2-7}) is
\begin{align}\label{2-8}
\frac{1}{m_2-1}\zeta(m_1+m_2-1)+O(|\varepsilon_2|)
+O(|\varepsilon_1+\varepsilon_2|).
\end{align}

Next we consider the second term. Since $m_1\leq 0$, $m_2\leq 0$, 
we see that $M_2({\mathbf m})=1-m_1-m_2$. Noting
$$
\binom{-m_2-\varepsilon_2}{k_2}=\frac{(-1)^{k_2}}{k_2!}
m_2(m_2+1)\cdots(m_2+k_2-1)+O(|\varepsilon_2|),
$$
when $k_2\neq 1-m_1-m_2$ we have
\begin{align}\label{2-10}
&\binom{-m_2-\varepsilon_2}{k_2}\zeta(-k_2)
\zeta(m_1+m_2+k_2+\varepsilon_1+\varepsilon_2)\\
&\;=\frac{(-1)^{k_2}}{k_2!}m_2(m_2+1)\cdots(m_2+k_2-1)
\zeta(-k_2)\zeta(m_1+m_2+k_2)\notag\\
&\;+O(|\varepsilon_2|)+O(|\varepsilon_1+\varepsilon_2|).\notag
\end{align}
When $k_2=1-m_1-m_2$, we have
$$
\zeta(m_1+m_2+k_2+\varepsilon_1+\varepsilon_2)=\frac{1}{\varepsilon_1+\varepsilon_2}
+\gamma+O(|\varepsilon_1+\varepsilon_2|),            
$$
and hence
\begin{align}\label{2-11}
&\binom{-m_2-\varepsilon_2}{k_2}\zeta(-k_2)
\zeta(m_1+m_2+k_2+\varepsilon_1+\varepsilon_2)\\
&\;=\binom{-m_2-\varepsilon_2}{1-m_1-m_2}\zeta(m_1+m_2-1)
\left(\frac{1}{\varepsilon_1+\varepsilon_2}
+\gamma+O(|\varepsilon_1+\varepsilon_2|)\right)\notag\\
&\;=\frac{(-1)^{1-m_1-m_2}}{(1-m_1-m_2)!}\frac{(m_2+\varepsilon_2)
(m_2+1+\varepsilon_2)\cdots(-m_1+\varepsilon_2)}{\varepsilon_1+\varepsilon_2}
\zeta(m_1+m_2-1)\notag\\
&\quad+\frac{(-1)^{1-m_1-m_2}}{(1-m_1-m_2)!}m_2(m_2+1)\cdots(-m_1)
\zeta(m_1+m_2-1)\gamma\notag\\
&\quad+O(|\varepsilon_2|)+O(|\varepsilon_1+\varepsilon_2|).\notag
\end{align}

Therefore from
\eqref{2-8}, \eqref{2-10} and \eqref{2-11} we obtain
\begin{align*}
\lefteqn{\zeta_2(m_1+\varepsilon_1,m_2+\varepsilon_2)}\\
&=\frac{1}{m_2-1}\zeta(m_1+m_2-1)\notag\\
&\;+\sum_{k_2=0}^{-m_1-m_2}
\frac{(-1)^{k_2}}{k_2!}m_2(m_2+1)\cdots(m_2+k_2-1)\zeta(-k_2)\zeta(m_1+m_2+k_2)
\notag\\
&\;+\frac{(-1)^{1-m_1-m_2}}{(1-m_1-m_2)!}\frac{(m_2+\varepsilon_2)                      
(m_2+1+\varepsilon_2)\cdots(-m_1+\varepsilon_2)}{\varepsilon_1+\varepsilon_2}           
\zeta(m_1+m_2-1)\notag\\                                                                
&\;+\frac{(-1)^{1-m_1-m_2}}{(1-m_1-m_2)!}m_2(m_2+1)\cdots(-m_1)                      
\zeta(m_1+m_2-1)\gamma\notag\\                                                          
&\;+O(|\varepsilon_2|)+O(|\varepsilon_1+\varepsilon_2|).\notag
\end{align*}
However when $k_2>-m_2$, the summands in the second term are $0$, 
since every summand has the factor $m_2(m_2+1)\cdots(m_2+k_2-1)=0$. Furthermore, 
the fourth term is also $0$, since $m_2(m_2+1)\cdots(-m_1)=0$. 
Therefore the assertion follows.
{\hfill $\square$}\\

At the point $(m_1,m_2)\in\mathbb{Z}_{\leq0}\times\mathbb{Z}$, we can also deduce the result  similar to Theorem \ref{Main4} by using the harmonic product \eqref{harmonic}. Hence, in the case $r=2$, we can give the Laurent series expansion at all points $\mathbf{m}\in\mathbb{Z}^2\setminus(\mathbb{Z}_{\geq1})^2$ with $O$-terms such as $O(|s_k-m_k|)$ or $O(|s_1+s_2-m_1-m_2|)$.

In the case $r=3$, there happens the situation when we cannot evaluate the limit value, if we use the above formula including $O$-terms for the double zeta values. For example, for $m_3\in\mathbb{Z}_{\leq0}$, we have
\begin{align}
\lefteqn{\zeta_3(s_1,s_2,s_3)}\notag\\
&=\frac{1}{s_3-1}\zeta_{2}(s_1,s_{2}+s_3-1)\notag\\
&\;+\sum_{k_3=0}^{M_3(\mathbf{m})}\binom{-s_3}{k_3}\zeta_{2}
    (s_1,s_{2}+s_3+k_3)\zeta(-k_3)\notag\\
&\;+O(|s_3-m_3|).\notag
\end{align}
By using the Laurent series expansion for $\zeta_2(s_1,s_2)$ with  $O$-terms, we can give the Laurent series expansion for $\zeta_3$ at the point $(m_1,m_2,m_3)\in\mathbb{Z}^2\times\mathbb{Z}_{\leq0}$. When $m_2\leq0$ and $m_3>0$, by using harmonic product
\begin{align}\label{harmonic3}
&\zeta_3(s_1,s_2,s_3)=\zeta(s_3)\zeta_2(s_1,s_2)-\zeta_2(s_1+s_3,s_2)\\
&\;\;-\zeta_3(s_1,s_3,s_2)
-\zeta_2(s_1,s_2+s_3)-\zeta_3(s_3,s_1,s_2),\notag
\end{align}
we can give the Laurent series expansion except for $m_3=1$. In the case $m_3=1$, the first term of the right-hand side of \eqref{harmonic3} contains the term $O(|s_2-m_2|/(s_3-1))$. Hence, in this case, we can not give the limit value. Therefore, generally we can not ignore the integral terms $I(s_1,\ldots,s_r;M_r(\mathbf{m})+1-\eta)$ in \eqref{LSE3} in order to get the limit value.


\end{document}